\documentclass[reqno]{amsart}
\usepackage{amssymb}
\usepackage{hyperref}

\title[ Riesz decomposition relative to $C^{\star}$-algebra homomorphisms  ]
{   Riesz decomposition relative to $C^{\star}$-algebra homomorphisms  }

\author[ A. Tajmouati, A. El Bakkali  and  S. Alaoui Chrifi ]
{  A. Tajmouati, A. El Bakkali  and  S. Alaoui Chrifi }

\address{ A. Tajmouati, S. Alaoui Chrifi\, \newline
Sidi Mohamed Ben Abdellah
 University,
 Faculty of Sciences Dhar Al Mahraz, Laboratory of Mathematical Analysis and Applications, Fez, Morocco.}
\email{abdelaziz.tajmouati@usmba.ac.ma}
\email{safae.alaouichrifi@usmba.ac.ma}

\address{ A. El Bakkali. \newline
Department of Mathematics.
University Chouaib Doukkali.
Faculty of Sciences Eljadida.
24000 Eljadida. Morocco.}
\email{aba0101q@yahoo.fr}

\subjclass[2010]{47A10, 47B06, 47C15}
\keywords{Homomorphism of Banach algebras. Riesz element in Banach algebra. Generalized Riesz element. Spectral projection. West-decomposition.}

\newtheorem{theorem}{Theorem}[section]
\newtheorem{definition}{Definition}[section]
\newtheorem{remark}{Remarks}
\newtheorem{lemma}{Lemma}[section]

\newtheorem{corollary}{Corollary}[section]
\newtheorem{example}{Example}
\newcommand{\pr} {{\bf   Proof: \hspace{0.3cm}}}
\begin{document}
\maketitle
\begin{abstract}
The aim of the present work is to study the Riesz decomposition relative to a $C^{\star}$-algebra homomorphism $T:\mathcal{A}\rightarrow \mathcal{B}$. We prove that under some conditions on $T$, $T$-Riesz elements can be decomposed into the sum of almost $T$-null element and quasi-nilpotent element. Also the so-called polynomial $T$-Riesz and generalized T-Riesz decompositions will be discussed.\\
\end{abstract}
\section{Introduction and Preliminaries }
Atkinson's and Ruston's characterizations in classical Fredholm theory have motivated several authors to introduce a general Fredholm theory in Banach algebra relative to a closed inessential ideal, this theory has been used by Smyth to generalize the well-known West decomposition in a $C^{\star}$-algebra, for more details of these decompositions we refer the reader to \cite{3, 4, 11, 12} as well as some characterizations of polynomially Riesz operators are given by authors in \cite{5,6,13}.\\

Let $\mathcal{A}$, $\mathcal{B}$ two complex Banach algebras with identity $e\neq0$, $T:\mathcal{A}\rightarrow \mathcal{B}$ be a homomorphism, $\mathcal{A}^{-1}$ the set of invertible elements in $\mathcal{A}$, $QN(\mathcal{A})$ the set of quasinilpotent elements in $\mathcal{A}$, $T^{-1}(0)$ the null space of the homomorphism $T$ and $Rad(\mathcal{B})$ the Jacobson radical of $\mathcal{B}$. It is clair that $T \mathcal{A}^{-1}\subset \mathcal{B}^{-1}$ and $T QN(\mathcal{A})\subset QN(\mathcal{A})$.\\

For any $a\in \mathcal{A}$ denotes $\sigma_{\mathcal{A}}(a)$, $\sigma_{\mathcal{B}}(Ta)$, $r_{\mathcal{A}}(a)$, $iso\sigma_{\mathcal{A}}(a)$ and $acc(\sigma_{\mathcal{A}}(a))=\sigma_{\mathcal{A}}(a)\backslash iso(\sigma_{\mathcal{A}}(a))$, the spectrum of a in $\mathcal{A}$, the spectrum of $Ta$ in $\mathcal{B}$, the spectral radius of a in $\mathcal{A}$, the set of isolated points and the accumulation points of $\sigma_{\mathcal{A}}(a)$ respectively. Recall that an ideal $\mathcal{J}$ is inessential ideal if for every $x\in \mathcal{J}$ the spectrum $\sigma_{\mathcal{A}}(a)$ is either finite or is a sequence converging to zero.\\

From Harte \cite{8} the set of $T$-Fredholm elements, $T$-Weyl elements, $T$-Browder elements and almost invertible $T$-Fredholm elements are define respectively by:
\begin{align*}
& T^{-1}\mathcal{B}^{-1}=\{a\in \mathcal{A}: Ta\in \mathcal{B}^{-1}\},\\
& T^{-1}(0)+\mathcal{A}^{-1}=\{c+d:  c\in T^{-1}(0) \>and\> d\in \mathcal{A}^{-1}\},\\
& T^{-1}(0)+_{com}\mathcal{A}^{-1}=\{c+d: c\in T^{-1}(0),\> d\in \mathcal{A}^{-1}\>and\> cd=dc\},\\
&and\> \{a\in \mathcal{A}: a\in T^{-1}\mathcal{B}^{-1} \>and\> 0\notin acc\sigma_{\mathcal{A}}(a)\}.
\end{align*}
Clearly we have
\begin{equation*}
Invertible \Rightarrow almost \text{     }invertible \text{     } T-Fredholm\> \Rightarrow\>T-Browder\> \Rightarrow\>T-Weyl\> \Rightarrow\>T-Fredholm.
\end{equation*}

For every $a\in \mathcal{A}$, let us define $T$-Fredholm spectrum, the $T$-Weyl spectrum, the $T$-Browder spectrum and the almost invertible $T$-Fredholm spectrum as follows respectively by:
\begin{align*}
 \sigma_{T}(a)                      &=\{\lambda \in \mathbb{C}: \lambda e-a\> is \>not\,\, T-Fredholm\}=\sigma_{\mathcal{B}}(Ta),\\
 \omega_{T}(a)                      &=\{\lambda \in \mathbb{C}: \lambda e-a\> is \>not\,\, T-Weyl\},\\
 \beta_{T}(a)                       &=\{\lambda \in \mathbb{C}:  \lambda e-a\> is \>not\,\, T-Browder\},\\
 and\>\>\> \sigma_{T}(a)\cup acc\sigma_{A}(a) &= \{\lambda \in \mathbb{C}:  \lambda e-a\> is \>not\,\,almost \>invertible \>T-Fredholm\}.
\end{align*}
Evidently
 $$\sigma_{T}(a)\subseteq \omega_{T}(a) \subseteq \beta_{T}(a) \subseteq \sigma_{T}(a)\cup acc\sigma_{\mathcal{A}}(a)\subseteq \sigma_{\mathcal{A}}(a).$$

On the other hand, let $K\subseteq \mathbb{C}$ be a compact subset, $\partial K$ the topological boundary of $K$ and $\eta K$ the connected hull of $K$ with $\mathbb{C}\backslash \eta K$ is the unique unbounded component of $\mathbb{C} \backslash K$, a hole of K is  a component of $\eta K\backslash K$,  thus $\eta K$ is the union of K and its holes.\\
It is known \cite{7} that for compact subsets $H,K\subseteq \mathbb{C}$
 $$\partial H\subseteq K\subseteq H \Longrightarrow \partial H\subseteq \partial K \subseteq K \subseteq H \subseteq \eta K= \eta H.$$
 If $K\subseteq \mathbb{C}$ is finite then $\eta K= K$. Consequently for compact subsets $H, K\subseteq \mathbb{C}$
 $$\eta K= \eta H \Longrightarrow H finite \Leftrightarrow K finite.$$
The homomorphism $T$ is said to has the Riesz property if for every $a\in \mathcal{A}$ we have
 $$a\in T^{-1}(0) \Longrightarrow acc(\sigma_{\mathcal{A}}(a))\subseteq \{0\}.$$
 Equivalently, $T$ satisfies the Riesz property if and only if $T^{-1}(0)$ is inessential ideal.\\
 And $T$ has the strong Riesz property, if $$\partial \sigma_{\mathcal{A}}(a) \subseteq \sigma_{T}(a) \cup iso(\sigma_{\mathcal{A}}(a)),$$
 Therefore $T$ satisfies the strong Riesz property if and only if $\sigma_{\mathcal{A}}(a) \subseteq \eta\sigma_{T}(a) \cup iso(\sigma_{\mathcal{A}}(a))$.\\
Recall that if $T$ is a bounded homomorphism and has the Riesz property then almost invertible $T$-Fredholm are $T$-Browder, in this case for an arbitrary $a\in \mathcal{A}$:
$$\beta_{T}(a) = \sigma_{T}(a)\cup acc(\sigma_{\mathcal{A}}(a)).$$
 According to \cite{10,14}, $a\in \mathcal{A}$ is $T$-Riesz element if $\lambda e-a$ is $T$-Fredholm for every $\lambda\in \mathbb{C}\backslash \{0\}$. We have the following characterization of $T$-Riesz elements:
 $$a \> is \>T-Riesz \>\> \Longleftrightarrow \>\>Ta\>is\>quasi-nilpotent \> in \> \mathcal{B}.$$
 The set of all $T$-Riesz elements is denoted by $T^{-1}QN(\mathcal{B})$.\\
 Now, let $T$ be a bounded homomorphism which has the Riesz property and \\ $\pi : \mathcal{A}\rightarrow \mathcal{A}/ T^{-1}(0) $ the quotient map.\\
 Write $\mathcal{R}_{T}(\mathcal{A})=\{a\in \mathcal{A}: \> \pi(a)\,\,\, is\,\,\, quasi-nilpotent\,\,\, in\,\,\, \mathcal{A}/ T^{-1}(0)\}$ the set of Riesz elements relative to $T^{-1}(0)$. It is easy to show that $$\mathcal{R}_{T}(\mathcal{A})\subseteq T^{-1} QN(\mathcal{B})$$
Furthermore if $a\in T^{-1} QN(\mathcal{B})$ then $acc(\sigma_{\mathcal{A}}(a))\subseteq \{0\}$.\\

This work focus on Fredholm theory relative to Banach algebra homomorphism introduced by Harte \cite{8, 9, 10, 14} in the special case of $C^{\star}$-algebras to extend West and Smyth decompositions also the decomposition of polynomially Riesz operators will be generalized .\\

\section{Main results}
 Throughout the rest of this paper $\mathcal{A}$ and $\mathcal{B}$ are two $C^{\star}$-algebras however proofs are due to Gelfand-Naimark's theorem [\cite{1}, theorem 6.2.20]  which states that every $C^{\star}$-algebra is isometrically $\star$-isomorph to a closed  self-adjoint subalgebra of $\mathcal{B}(\mathcal{H})$ for some Hilbert space $\mathcal{H}$.

 \begin{definition} \cite{14}
 Let $T:\mathcal{A}\longrightarrow \mathcal{B}$ be a homomorphism, an element $a\in \mathcal{A}$ is said $T$-null if $Ta=0$ and almost $T$-null if $Ta\in Rad(\mathcal{B})$.
 \end{definition}

Because $C^{\star}$-algebras are semi-simple then $Rad(\mathcal{B})=0$, in this case $T$-null and almost $T$-null elements are the same.\\

Now, we are in position to establish the following theorem, the technique used for dealing with the proof of [\cite{3}, $C^{\star}$.2.3 Theorem] can also be used for showing the following theorem.
 \begin{theorem}$\label{2.1}$
 Let $T:\mathcal{A}\rightarrow \mathcal{B}$ be a bounded homomorphism which has Riesz property. Then every $a\in T^{-1} QN(\mathcal{B})$ has the decomposition $a=c+d$ with d is a normal element of $T^{-1}(0)$ and $c\in QN(\mathcal{A})$, i.e:
 $$T^{-1} QN(\mathcal{B})= T^{-1}(0)+ QN(\mathcal{A}).$$
 \end{theorem}
\pr\\
  It is clair that $T^{-1}(0)+ QN(\mathcal{A})\subseteq T^{-1} QN(\mathcal{B})$ so it remains to show the other inclusion, in fact let $a\in T^{-1} QN(B)$,
  if $\sigma_{\mathcal{A}}(a)=\{0\}$ or $\sigma_{\mathcal{A}}(a)$ is finite then the result is trivial. To avoid triviality let assume that $\sigma_{\mathcal{A}}(a)$ has an infinite number of isolated points, set $\sigma_{\mathcal{A}}(a)=(\lambda_{k})_{k\geq1}\cup\{0\}$ with $\{|\lambda_{k}|\}_{1}^{\infty}$ is a decreasing sequence of non zero isolated points of $\sigma _{\mathcal{A}}(a)$ such that $\lim \lambda_{k}=0$.\\
  Consider $P_{j}$ the spectral projection of $a$ associated to the spectral set $\{\lambda_{j}\}$ given by:
  $$P_{j} = \frac{1}{2\pi i}\int_{\Gamma_{j}}(\lambda e_{\mathcal{A}}-a)^{-1}d\lambda$$
  where $\Gamma_{j}$ is a circle of center $\lambda_{j}$ which contains no other points of $\sigma_{\mathcal{A}}(a)$. By [\cite{10}, theorem 1.1] we have
  $$TP_{j}=\frac{1}{2\pi i}\int_{\Gamma_{j}}(\lambda e_{\mathcal{A}}-Ta)^{-1}d\lambda= 0 (\mbox{ since } \,\,\,\,\lambda_{j}\notin \sigma_{\mathcal{B}}(Ta)\>\> for\>every\> j=1,2,...)$$
  Therefore $P_{j}\in T^{-1}(0)$ for every j=1,2,... ,
  set $S_{n}=\sum_{k=1}^{n}{P_{k}}$ and let $q_{n}$ be a self-adjoint projection in $T^{-1}(0)$ such that $q_{n}(\mathcal{H})=S_{n}(\mathcal{H})$.\\
  Write $E_{i}=q_{i}-q_{i-1}$ with $q_{0}=0$ and set $d_{n}=\sum_{k=1}^{n}{ \lambda_{k} E_{k}}$, then $d_{n}\in T^{-1}(0)$. Since $\{E_{i}\}_{1}^{n}$ is a sequence of self-adjoint projections then $d_{n}$ is a normal element of $T^{-1}(0)$ for every n and $d_{n}\longrightarrow d\in T^{-1}(0)$ when $n\rightarrow \infty$.\\

  Considering the decomposition $\mathcal{H}=q_{n}(\mathcal{H})\bigoplus (1-q_{n})(\mathcal{H})$, the matrix formes associated to this decomposition are:\\
  $d_{n}= \left(\begin{array}{c c c}
  k_{n}\>\>\> 0\\
               \\
  0\>\>\>\>\>0
  \end{array} \right)$
with $\sigma(d_{n})=\sigma(k_{n})=\{\lambda_{k}\}_{1}^{n}$\\
and
  $a= \left(\begin{array}{c c c}
  k_{n}\>\>\> *\\
               \\
  0\>\>\>\>\>Z
  \end{array} \right)$
  with
  $\sigma(a)=\sigma(k_{n})\cup \sigma(Z)$, consequently $\sigma(Z)=\sigma(a)\setminus \{\lambda_{k}\}_{1}^{n}.$\\
Now, it suffice to show that $c=a-d$ is quasi-nilpotent, to show this let $\lambda\in \mathbb{C}\backslash \{0\}$. Without lose generality we may assume that $\lambda \notin  \sigma(a-d_{n})$, then $\lambda e-a+d_{n}$ is invertible.\\
We have  $a-d_{n}= \left(\begin{array}{c c c}
  0\>\>\>\>\> *\\
               \\
  0\>\>\>\>\>Z
  \end{array} \right)$
  then $\sigma(a-d_{n})=\sigma(Z)=\sigma(a)\setminus \{\lambda_{k}\}_{1}^{n}$\\
  Observe that $(\lambda e-a+d_{n})^{-1}=\left(\begin{array}{c c c}
  \frac{1}{\lambda}I_{n}\>\>\>\>\>\>\>\>\>\>\>\>*\>\>\>\>\>\>\>\> \\
               \\
  \>\>0\>\>\>\>\>\>\>(\lambda-Z)^{-1}
  \end{array} \right)$
     and $d-d_{n}=\left(\begin{array}{c c c}
  0\>\>\>\>\> 0\\
               \\
  0\>\>\>\>\>h_{n}
  \end{array} \right)$ \\
  where $I_{n}$ is the identity matrix and $\sigma(h_{n})=\sigma(d)\setminus \{\lambda_{k}\}_{1}^{n}$.\\
Thus $(\lambda-a+d_{n})^{-1}(d-d_{n})=\left(\begin{array}{c c c}
  0\>\>\>\>\>\>\>\>\>\>\>\>\>\>\>\>\>\>\>\>*\>\>\>\>\>\>\> \\
               \\
  \>\>0\>\>\>\>\>\>\>h_{n}(\lambda-Z)^{-1}
  \end{array} \right).$\\
 Assume that $\|h_{n}\|< \|(\lambda-Z)^{-1}\|^{-1}$ then\\
\begin{align*}
r_{\mathcal{A}}((\lambda e-a+d_{n})^{-1}(d-d_{n})) &\leq r_{\mathcal{A}}(h_{n}(\lambda-Z)^{-1})\\
&\leq \|h_{n}\| \|(\lambda-Z)^{-1}\|\\
&< 1.
\end{align*}
Therefore $1+(\lambda e-a+d_{n})^{-1}(d-d_{n})$ is invertible, which means that $\lambda e-a+d$ is invertible and hence $\lambda \notin \sigma(a-d)$, consequently $\sigma(c)={0}$. \,\,\,\,\,\,\,\,$\blacksquare$

\begin{remark}
If in addition $T$ is onto then $\mathcal{R}_{T}(\mathcal{A})=T^{-1}QN(\mathcal{B})$, in this case the decomposition in theorem \ref{2.1} follows immediately from \cite{11}.
\end{remark}

\begin{example}
Let $\mathcal{H}$ be a Hilbert space, $\mathcal{B}(\mathcal{H})$ the algebra of bounded linear operators on $\mathcal{H}$, $\mathcal{K}(\mathcal{H})$, $\mathcal{QN}(\mathcal{H})$ and $\mathcal{R}(\mathcal{H})$ the closed ideal of compact operators, the class of quasi-nilpotent operators and the class of Riesz operators respectively. Consider $\mathcal{C}(\mathcal{H})=\mathcal{B}(\mathcal{H})/\mathcal{K}(\mathcal{H})$ the Calkin algebra and $\pi : \mathcal{B}(\mathcal{H})\longrightarrow \mathcal{C}(\mathcal{H})$ the canonical bounded homomorphism.\\
It is clear that $\pi$ has the Riesz property since $\pi^{-1}(0)=\mathcal{K}(\mathcal{H})$ is an inessential ideal.\\
In the other hand if $T\in \mathcal{B}(\mathcal{H})$, the
Atkinson's characterization said that $T$ is Fredholm if and only if $\pi(T)$ is invertible in the Calkin algebra, thus
$$\sigma_{\mathcal{C}(\mathcal{H})}(\pi(T))=\sigma_{\pi}(T)=\{\lambda\in\mathbb{C}: T-\lambda e \>\> is\>not\> Fredholm\}.$$
According to Ruston characterization, $T$ is Riesz operator if and only if  $\pi(T)$ is quasi-nilpotent in the Calkin algebra, so by West \cite{12} we have the following decomposition
$$\mathcal{R}(\mathcal{H})=\pi^{-1}(0)+\mathcal{QN}(\mathcal{H}).$$
\end{example}

\begin{definition}\cite{14}
Let $a\in \mathcal{A}$, $a$ is called polynomial $T$-Riesz element if there exist a non-zero polynomial P such that $P(a)$ is $T$-Riesz.\\
Note $poly^{-1}T^{-1}(QN(\mathcal{B}))$ the set of all polynomial $T$-Riesz elements.
\end{definition}

\begin{lemma}$\label{lemma 2.1}$
Let $a\in poly^{-1}T^{-1}(QN(\mathcal{B}))$ with $T$ is a bounded homomorphism which has the strong Riesz property then
\begin{equation}
\sigma_{T}(a)=\omega_{T}(a)=\beta_{T}(a)=P^{-1}(0),
\end{equation}
For a polynomial P of minimal degree, and $acc(\sigma_{\mathcal{A}}(a))\subseteq \omega_{T}(a)$.
\end{lemma}

\pr\\
 By [\cite{14}, theorem 11.1] (1) is hold, \\
 In this case we have $\omega_{T}(a)=\sigma_{T}(a)\cup acc\sigma_{\mathcal{A}}(a)$ then obviously $acc(\sigma_{\mathcal{A}}(a))\subseteq \omega_{T}(a)$.\\

The following theorem gives the decomposition of polynomially $T$-Riesz element in $C^{\star}$-algebra.
 \begin{theorem}
 Let T be a bounded homomorphism which has the strong Riesz property then $$a\in poly^{-1}T^{-1}(QN(\mathcal{B}))\>\>\>\Longleftrightarrow \>\>\> \omega_{T}(a)\>is\>finite$$
In this case the following decomposition $a=c+d+f$ is hold with
  $c\in T^{-1}(0)$,
  $d\in QN(\mathcal{B})$
  and $\sigma(f)=\omega_{T}(a).$
 \end{theorem}

\pr\\
 Assume that $\omega_{T}(a)$ is finite, claim that $acc\sigma_{\mathcal{A}}(a)\subseteq \omega_{T}(a)$.\\
 In fact, let $\lambda\in \sigma_{\mathcal{A}}(a)\backslash \omega_{T}(a)$, then $\lambda e-a $ is $T$-Fredholm but not invertible, thus $\lambda\in \partial \sigma_{\mathcal{A}}(a)$, if it is not then $\lambda\in int(\sigma_{\mathcal{A}}(a))$, in this case $\sigma_{T}(a)$ would be infinite which is impossible.\\
 So $\lambda\in \partial \sigma_{\mathcal{A}}(a)\backslash \sigma_{T}(a)\subseteq iso\sigma_{\mathcal{A}}(a)$, and therefore
 $$acc(\sigma_{\mathcal{A}}(a))\subseteq \omega_{T}(a).$$
 Set $\sigma_{T}(a)=\omega_{T}(a)=\{\lambda_{1},\lambda_{2},...\lambda_{n}\}$, so $acc\sigma_{\mathcal{A}}(a)\subseteq \{\lambda_{1},\lambda_{2},...\lambda_{n}\}$.\\
 Let now define a collection $\{\triangle_{1},\triangle_{2},...\triangle_{n}\}$ of subsets of $\sigma_{\mathcal{A}}(a)$ in the following way:\\
 $(i)\>\>\>\bigcup_{i=1}^{n}\triangle_{i}=\sigma_{\mathcal{A}}(a),$\\
 $(ii)\>\>\> \triangle_{i}\cap \triangle_{j}=\emptyset \>\>\>if\>i\neq j,$\\
 $(iii)\>\>\>\lambda_{i}\in \triangle_{i} \>\>\> for \>\>i=1,2,...,n,$ \\
 For every i=1,2,..,n let $N_{i}$ be a neighborhood of $\triangle_{i}$ which contains no other points of $\sigma_{\mathcal{A}}(a)$ than $\lambda_{i}$ and let $P_{i}$ the spectral projection of $a$ associated to $\triangle_{i}$ for i=1,2,...,n
 $$P_{i} = \frac{1}{2\pi i}\int_{\partial N_{i}}(\lambda e_{\mathcal{A}}-a)^{-1}d\lambda.$$
 We have the decomposition $\mathcal{H}=R(P_{1})\bigoplus R(P_{2})\bigoplus ...\bigoplus R(P_{n})$, thus we have
  $$a=\bigoplus_{i=1}^{n}(a_{i}+\lambda_{i}e_{i})\,\,\,\,\,\, \mbox{ and } \,\,\,\,\,\,Ta=\bigoplus_{i=1}^{n}(Ta_{i}+\lambda_{i}e_{i})$$
 Therefore $\sigma_{\mathcal{B}}(Ta)=\bigcup_{i=1}^{n}\sigma_{\mathcal{B}}(Ta_{i}+\lambda_{i}e_{i})$ and consequently $\sigma_{\mathcal{B}}(Ta_{i})=\sigma_{T}(a_{i})=\{0\}$ for every $i=1,2,...,n$.
 Then $a_{i}$ is $T$-Riesz for every $i=1,2,...,n$.\\
 Next, using theorem \ref{2.1} we have $a_{i}=c_{i}+d_{i}$ where $c_{i}\in QN(\mathcal{A})$ and $d_{i}\in T^{-1}(0)$.
 Consequently
  \begin{align*}
  a &= \bigoplus_{i=1}^{n}(c_{i}+d_{i}+\lambda_{i}e_{i})\\
    &=c+d+f
  \end{align*}
  where $c=\bigoplus_{i=1}^{n}c_{i}$, $d=\bigoplus_{i=1}^{n}d_{i}$ and $f=\bigoplus_{i=1}^{n}\lambda_{i}e_{i}$.
   so $c\in QN(\mathcal{A})$, $d\in T^{-1}(0)$ \\ and $\sigma(f)=\omega_{T}(a)$.\\
   Now, set $P(\lambda)\equiv \prod_{i=1}^{n}(\lambda-\lambda_{i})$ thus
   \begin{align*}
  P(a) &=\prod_{i=1}^{n}(a-\lambda_{i}e_{i})\\
    &=\prod_{i=1}^{n}( \bigoplus_{j=1}^{n}(c_{j}+d_{j}+\lambda_{j}e_{j})-\lambda_{i}e_{i})\\
    &=\bigoplus_{j=1}^{n}( \prod_{i=1}^{n}(c_{j}+d_{j}+\lambda_{j}e_{j}-\lambda_{i}e_{i}))\\
    &= L+ R
  \end{align*}
 With $L\in T^{-1}(0)$ and $R= c_{j}^{n}+ \mu_{n-1}c_{j}^{n-1}+...+\mu_{1}c_{j}$ for some $\mu_{i}\in \mathbb{C}$\,\,\,\,\,\, $i=1,2,...,n-1$.\\
 It is clear that $R\in QN(\mathcal{A})$
 and $$\sigma_{T}(P(a))=\sigma_{\mathcal{B}}(TP(a))=\sigma_{T}(R)={0}$$
 Hence $P(a)\in T^{-1}QN(\mathcal{B})$, thus the proof is complete.\,\,\,\,\,\, $\blacksquare$\\

Similarly as in [\cite{5}, definition 1.2] we have the following

\begin{definition}
Let $a\in \mathcal{A}$, $a$ is said generalized $T$-Riesz element if there exist a finite subset $\Omega\subseteq \mathbb{C}$ such that:
$$\sigma_{T}(a)=\omega_{T}(a)=\beta_{T}(a)=\Omega.$$
Write $GT^{-1}QN(\mathcal{B})$ the set of generalized $T$-Riesz elements.
\end{definition}

To conclude this section we have the following result and remark.
\begin{corollary}
Let $T:\mathcal{A}\longrightarrow \mathcal{B}$ be a bounded homomorphism which has the Riesz property then every $a\in GT^{-1}QN(\mathcal{B})$ has the decomposition $a=d+c$ with $d$ is normal element of $T^{-1}(0)$ and $\sigma(c)=\Omega$.
\end{corollary}

\begin{remark}
1)- Obviously $T^{-1}QN(B)\subseteq poly^{-1}T^{-1}QN(\mathcal{B})$ and if $a\in \mathcal{A}$ is polynomial $T$-Riesz for a polynomial $P$ which has no zero in $\sigma_{\mathcal{A}}(a)$, then a is $T$-Riesz.\\
2)- By lemma \ref{lemma 2.1} and theorem \ref{2.1} we have $$GT^{-1}QN(\mathcal{B})=poly^{-1}T^{-1}QN(\mathcal{B}).$$\\
\end{remark}


\begin{thebibliography}{99}
\bibitem{1} B. Aupetit, \textit{A Primer on Spectral Theory}, Springer-Verlag, Berlin (1991).\par

\bibitem{2} B. Aupetit, \textit{Proprietes spectrales des algebres de Banach}, in Vol. 735 ofLecture Notes in Mathematics,
Springer Verlag, Heidelberg (1979) .

\bibitem{3} B.A. Barnes, G.J.Murphy, M.R.F Smyth and T.T west, \textit{Riesz and Fredholm theory in Banach algebras}, London Pitman, 1982.

\bibitem{4} S.R. Caradus, W.E.Pfaffenberger and Bertran Yood, \textit{calkin algebra and algebras of operators on Banach spaces}, Marel Dekker INC, New york, 1974.\par

\bibitem{5} A. Jeribi and N. Moalla, \textit{Fredholm operators and Riesz theory for polynomially compact operators}, Acta Applicandae Mathematicae,
vol. 90, no. 3, pp. 227-247, 2006.

\bibitem{6} Y. M. Han, S. H. Lee, and W. Y. Lee, \textit{On the structure of polynomially compact operators}, Math. Z. 232,
257 (1999).\par

\bibitem{7} R.E Harte and A.W.Wicksted, \textit{Boundaries, Hulls and spectral mapping theorems}, Royal Irish Academy Vol. 81A, No. 2 (1981), pp. 201-208.\par

\bibitem{8} R.E Harte, \textit{Fredholm Theory Relative to a Banach Algebra Homomorphism}, Math. Z. 179, 431-6.\par

\bibitem{9} R.E Harte, \textit{Invertibility and singularity for bounded linear operators}, Marcel Dekker, New york 1988.\par

\bibitem{10} T. Mouton and H.Raubenheimer \textit{ More on Fredholm theory relative to Banach algebra homomorphisms } Math. Proc. R. Ir. Acad. 93A (1993) 17–25.\par

\bibitem{11} M.R.F Smyth, \textit{Riesz theory in Banach algebra}, Mathematsche Zeitschrift 145(2), 144(1975).\par

\bibitem{12} T.T West, \textit{The decomposition of Riesz operators}, Proc.London.Soc 6,737-752 (1966).\par

\bibitem{13} $\check{Z}ivkovi\acute{c}-Zlatanovi\acute{c}\> S\check{C}, Djordjevi\acute{c} \>DS, Harte\>RE, Duggal\>BP$, \textit{On polynomially Riesz
operators},Filomat. 2014;28:197–205.\par

\bibitem{14}  $S.\breve{C}. \breve{Z}ivkovi\acute{c}-Zlatanovi\acute{c}, D.S. Djordjevi\acute{c}, R.E. Harte$, \textit{Ruston, Riesz and perturbation classes}, J. Math. Anal. Appl. 389 (2012) 871–886.\par


\end{thebibliography}
\end{document}